# Partially Paradoxist Smarandache Geometries

**Howard Iseri**
Department of Mathematics and Computer Information Science
Mansfield University
Mansfield, PA 16933
**hiseri@mnsfld.edu**

Abstract: A paradoxist Smarandache geometry combines Euclidean, hyperbolic, and elliptic geometry into one space along with other non-Euclidean behaviors of lines that would seem to require a discrete space. A class of continuous spaces is presented here together with specific examples that exhibit almost all of these phenomena and suggest the prospect of a continuous paradoxist geometry.

**Introduction**

Euclid's parallel postulate can be formulated to say that given a line *l* and a point *P* not on *l*, there is exactly one line through *P* that is parallel to *l*. An axiom is said to be **Smarandachely denied**, if it, or one of its negations, holds in some instances and fails to hold in others within the same space. For example, Euclid's parallel postulate would be Smarandachely denied in a geometry that was both Euclidean and non-Euclidean, or non-Euclidean in at least two different ways. A **Smarandache geometry** is one that has at least one Smarandachely denied axiom, and a paradoxist Smarandache geometry, to be described later, denies Euclid's parallel postulate in a somewhat exhaustive way.

Euclid's parallel postulate does not hold in the standard non-Euclidean geometries, the hyperbolic geometry of Gauss, Lobachevski, and Bolyai and the elliptic geometry of Riemann. These are special cases of the two-dimensional manifolds of Riemannian geometry. Here the three types of geometry are characterized by the Gauss curvature, negative curvature for hyperbolic, zero curvature for Euclidean, and positive curvature for elliptic. In general, the curvature may vary within a particular Riemannian manifold, so it is possible that the geodesics, the straightest possible curves, will behave like the lines of Euclidean geometry in one region and like the lines of hyperbolic or elliptic geometry in another. We would expect, therefore, to find geometries among the Riemannian manifolds that Smarandachely deny Euclid's parallel postulate. The models presented here will suggest specific examples, but explicit descriptions would be far from trivial.

We will bypass the computational complexities of Riemannian manifolds by turning to a class of geometric spaces that we will call Smarandache manifolds or S-manifolds. S-manifolds are piecewise linear manifolds topologically, and they have geodesics that exhibit elliptic, hyperbolic, and Euclidean behavior similar to those in Riemannian geometry, but that are much easier to construct and describe.

The idea of an S-manifold is based on the *hyperbolic paper* described in **[2]** and credited to W. Thurston. There, the negative curvature of the hyperbolic plane is visualized by taping together seven triangles made of paper (see Figures 2a and 2b). Squeezing seven equilateral triangles around a vertex, instead of the usual six seen in a tiling of the plane, forces the paper into a flat saddle shape with the negative curvature concentrated at the center vertex. By utilizing these "curvature singularities," our S-manifolds can be flat (i.e., Euclidean) everywhere else.

**Smarandache manifolds**

A **Smarandache manifold** (or **S-manifold**) is a collection of equilateral triangular disks (triangles) where every edge is shared by exactly two triangles, and every vertex is shared by five, six, or seven triangles. The points of the manifold are those of the triangular disks, including all the interior points, edge points, and vertices. Lines (geodesics) in the manifold are those piecewise linear curves with the following properties. They are straight in the Euclidean sense within each triangular disk and pair of adjacent triangular disks (since two triangles will lie flat in the plane). Across a vertex, a line will make two equal angles (two 150º angles for five triangles, two 180º angles for six triangles, and two 210º angles for seven).

**Elliptic Vertices --- five triangles**

There are five equilateral triangles around an elliptic vertex in an S-manifold. We can take a region around an elliptic vertex and lay it flay by making a cut as in Figure 1a. Note that the lines are straight within any pair of adjacent triangles, although the lines appear to bend at the vertex and across the cut. This is only because we have made a cut and flattened the surface. In the paper model shown in Figure 1b, the lines curve, but only in a direction perpendicular to the surface. In other words, the lines are as straight as possible and bend only as they follow the surface. The two lines that do not pass through the central vertex pass through three adjacent triangles, which would lie flat in the plane, and so are straight in the Euclidean sense. Note that the fact that the third triangle is shared by both lines forces them to intersect. The middle line runs along an edge of a triangle and passes through an elliptic vertex, so it bisects the opposite triangle making two 150º angles (or two-and-a-half triangles). In general, lines passing on either side of an elliptic vertex will turn towards each other.

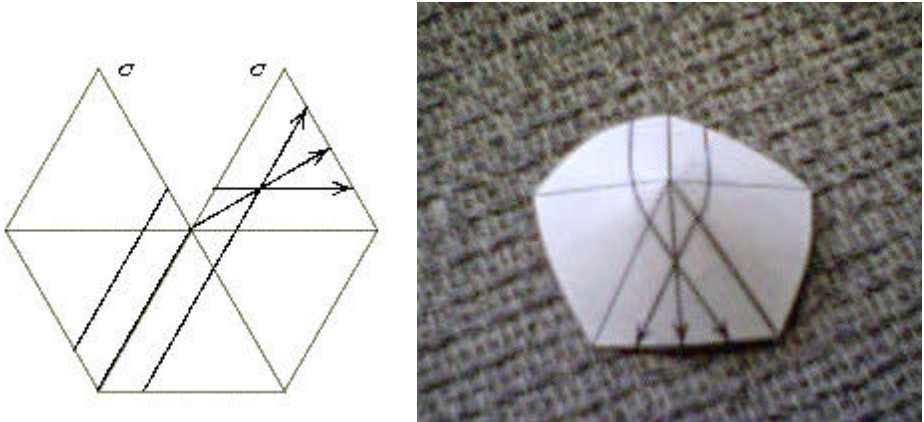

Figures 1a and 1b. Lines near an elliptic vertex.

**Hyperbolic Vertices --- seven triangles**
There are seven triangles around a hyperbolic vertex. We can lay a region around a hyperbolic vertex flat after making cuts as shown in Figure 2a. The middle line runs along an edge, so it bisects the opposite triangle (and has 210º, or three-and-a-half triangles, on either side of it). The two lines on either side pass through three adjacent triangles, and are straight as in the elliptic case. Note that the third triangles here are separated by another triangle, so lines passing on either side of a hyperbolic vertex turn away from each other.

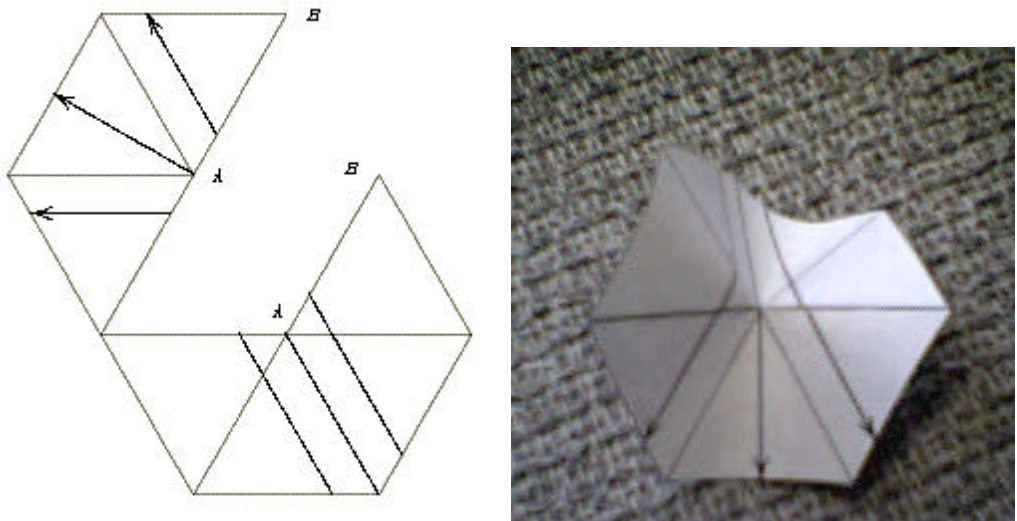

Figures 2a and 2b. Lines near a hyperbolic vertex.

**Paradoxist geometries**

We will say that a point *P* not on a line *l* is **Euclidean** with respect to *l*, if there is exactly one line through *P* that is parallel to *l*. *P* is **elliptic** with respect to *l*, if there are no parallels through *P*. If there are at least two parallels through *P*, then it is **hyperbolic**. Furthermore, if *P* is hyperbolic with respect to *l*, then it is **finitely hyperbolic**, if there are only finitely many parallels, and it is **regularly hyperbolic**, if there are infinitely many parallels and infinitely many non-parallels. Finally, if there are infinitely many parallels and only finitely many non-parallels, then *P* is **extremely hyperbolic**, and if all the lines through *P* are parallel, then *P* is **completely hyperbolic**.

Smarandache called a geometry **paradoxist** if there are points that are elliptic, Euclidean, finitely hyperbolic, regularly hyperbolic, and completely hyperbolic [1]. We will add extremely hyperbolic to the definition of a paradoxist geometry. We will also say that a geometry is **semi-paradoxist**, if it has Euclidean, elliptic, and regularly hyperbolic points, and if it lacks only finitely hyperbolic points we will call it **almost paradoxist.**

**A Semi-Paradoxist Model**

This model is constructed by taking a hyperbolic and an elliptic vertex adjacent to each other and surrounding them with Euclidean vertices to form a space that is topologically equivalent to the plane. A part of it is shown in Figures 3a and 3b. Let *l* be the line through *O*. With respect to *l*, we see that the point *P* is Euclidean. The line through *P* shown is parallel to *l*, and any other line through *P* clearly intersects *l*, since the region to the right and left is essentially Euclidean.

The point *Q* is elliptic with respect to *l*. The line shown intersects *l*, as would any other line through *Q*.

The point *R* is regularly hyperbolic. The lines shown are parallel to *l*, and these separate the other infinitely many parallels from the infinitely many non-parallels.

This S-manifold can be turned into a Riemannian manifold by smoothing the two curvature singularities. The lines shown in Figures 3a and 3b would stay the same, and only those geodesics passing near the singularities would be affected by the change.

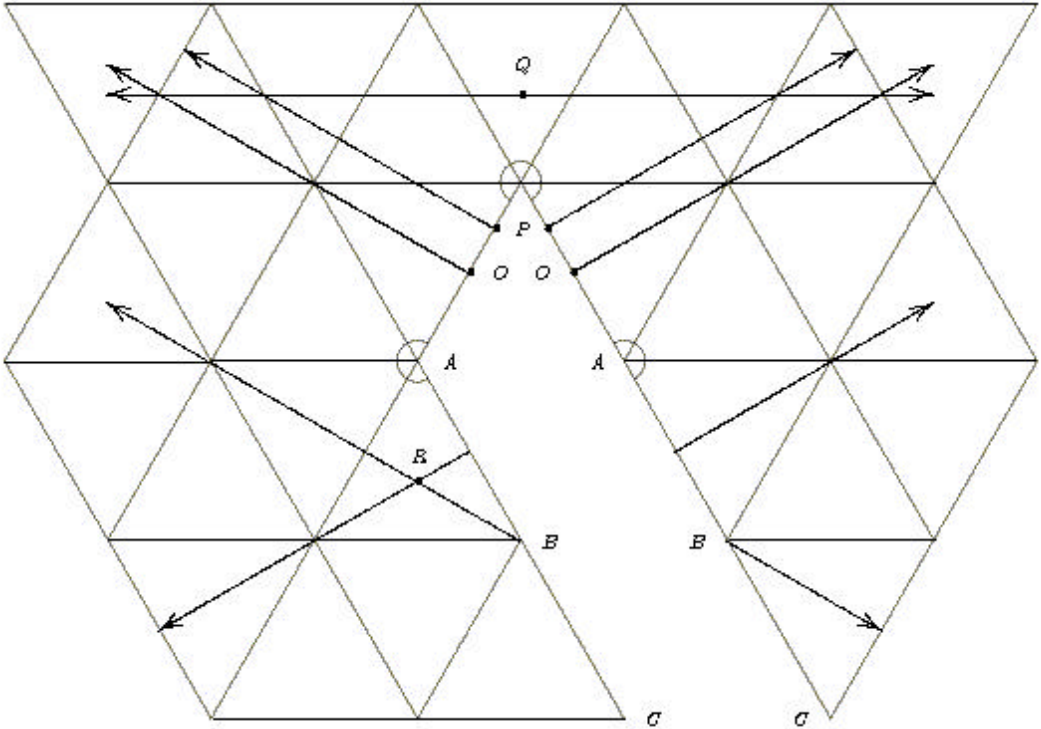

Figure 3a. Lines in the semi-paradoxist model.

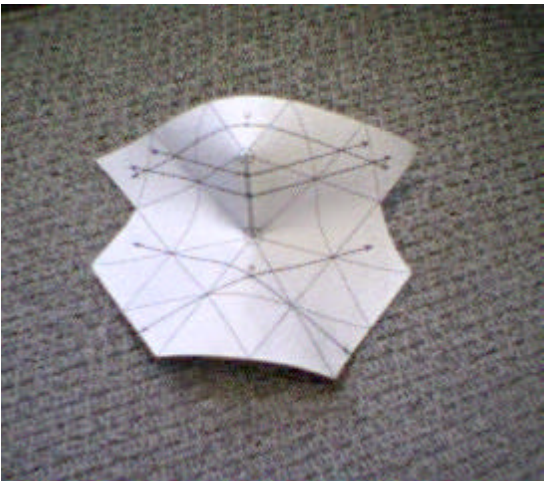

Figure 3b. Lines in the semi-paradoxist model.

**An Almost Paradoxist Model**
A greater variety in the types of hyperbolic points can be found in an S-manifold with more hyperbolic vertices. This model has at its center an elliptic vertex surrounded by five more elliptic vertices. Five Euclidean vertices then surround these elliptic vertices (see Figures 4a and 4b) to form a cylinder with a cone on top of it. We will call this the silo.

The line *l* runs around the cylinder (it is a circle). With respect to the line *l*, the point *P* is Euclidean, and the point *R* is elliptic.

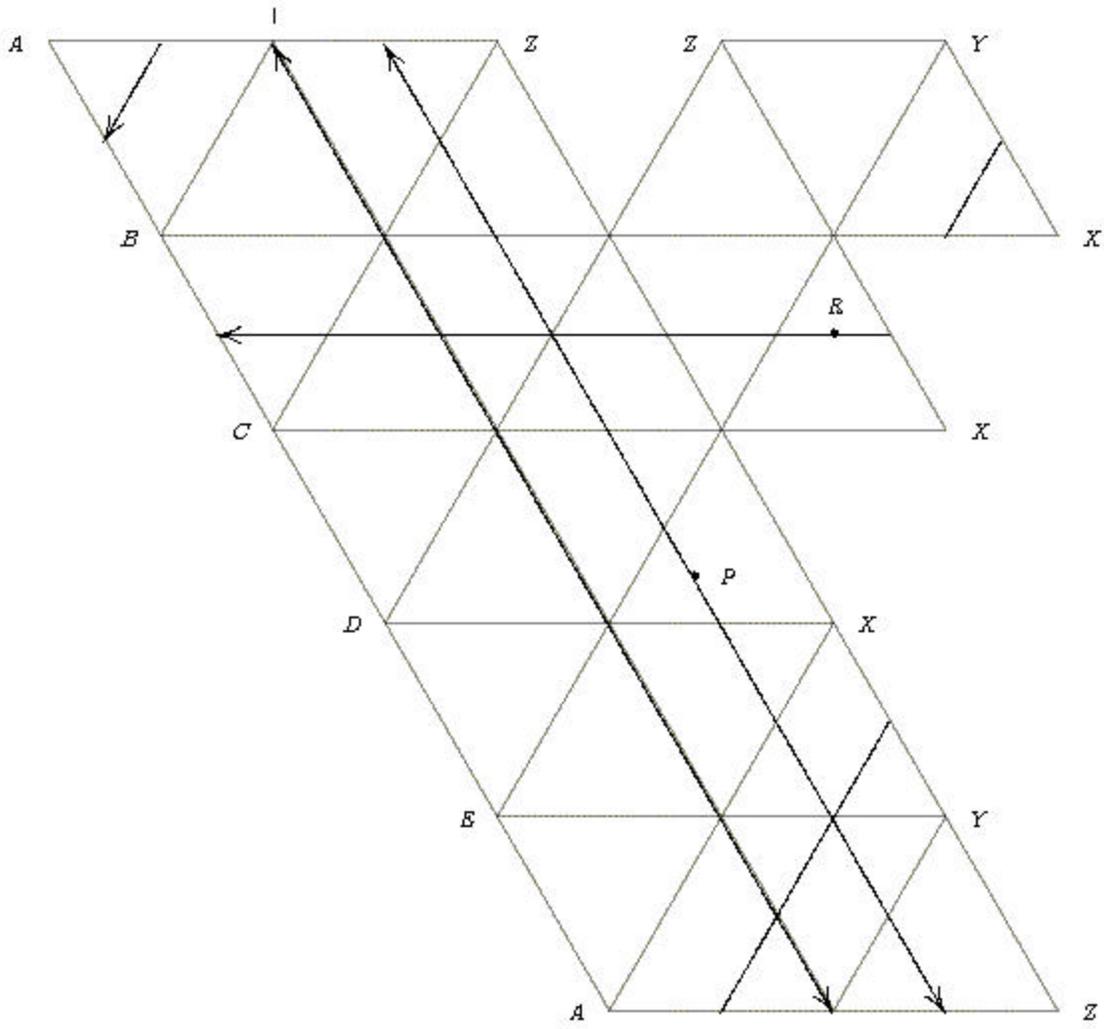

Figure 4a. Lines in the silo of the almost paradoxist model.

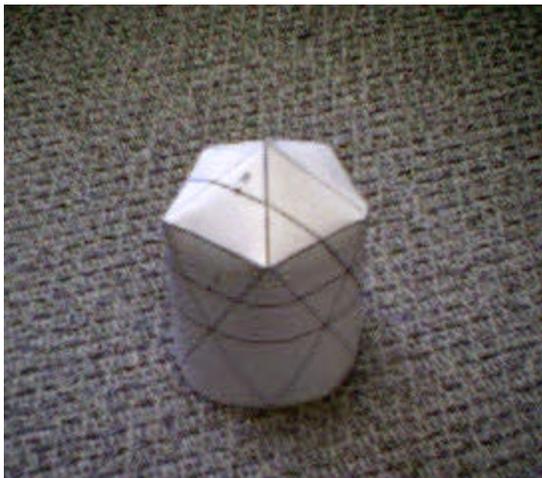
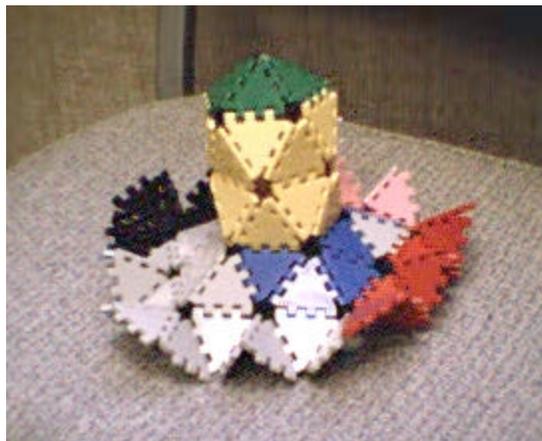

Figures 4b and 4c. Lines in the silo of the almost paradoxist model, and the hyperbolic region around silo.

The entire model is topologically equivalent to the plane, and is completed by extending the bottom of the silo with hyperbolic vertices. A model made of *ZAKS blocks* in Figure 4c shows some of the hyperbolic region extending from the bottom of the silo.

Examples of various types of hyperbolic points are shown in Figures 5a and 5b, which shows some of the hyperbolic region and the bottom of the silo. The line *l* mentioned previously is at the top. With respect to *l*, the point *Q* is regularly hyperbolic. The two lines shown are parallel to *l*, and they separate the parallels from the non-parallels. Out further into the hyperbolic region is the point *Q* . The line shown passing through *Q* and the vertex *I* intersects *l*. Any line through *Q* that misses the vertex *I* will lie outside of the two dotted lines, and these will miss the silo entirely. Since only one line through *Q* intersects *l*, it is an extremely hyperbolic point. The nearby point *Q* is completely hyperbolic. We can see this by noticing that the line through *Q* and *I* will follow the dotted line to the left and miss the silo. All the lines through *Q* to the left of this will also miss the silo. Any line to the right will miss the vertex *I*, and will run just to the right of the line through *Q* and *I* until it misses the vertex *F* and turns to the right. These lines will also miss the silo.

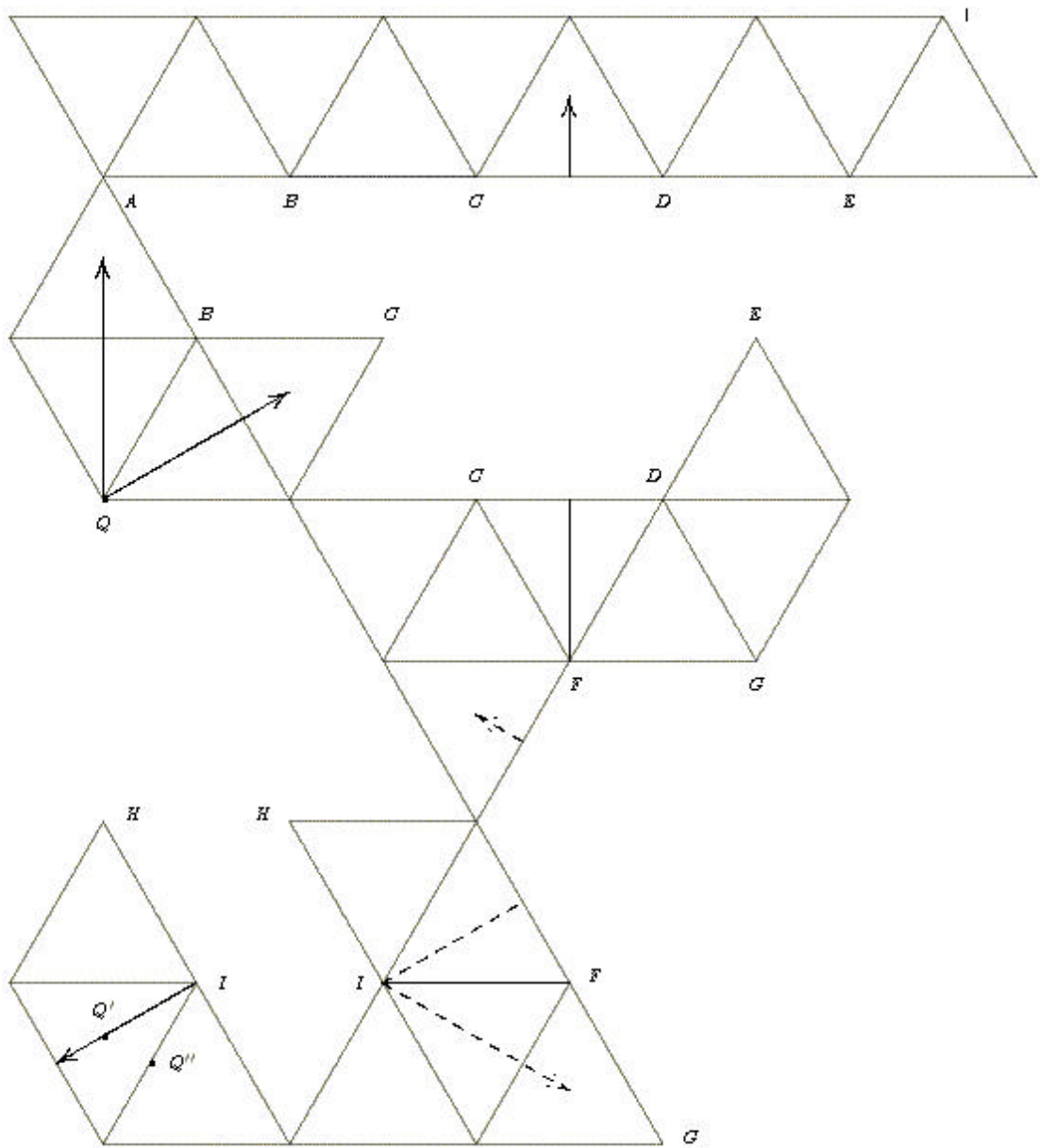

Figure 5a. Lines in the hyperbolic region near the silo in the almost paradoxist model.

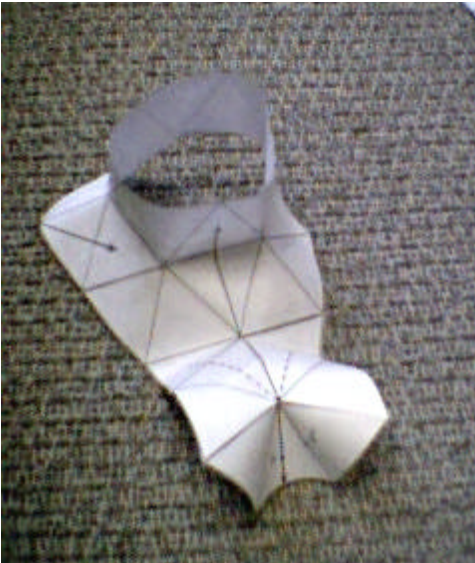
Figure 5b. Lines in the hyperbolic region near the silo in the almost paradoxist model.

Since this model has elliptic, Euclidean, regularly hyperbolic, extremely hyperbolic, and completely hyperbolic points, it is almost paradoxist. Note that it follows from the existence of extremely and completely hyperbolic points that there are pairs of points that do not lie on a single line. This model is connected, however, and there is always a finite sequence of line segments that connect any particular pair of points.

**Final Remarks**
It is relatively easy to construct an S-manifold that is almost paradoxist. The most interesting prospect, however, is the possibility of an S-manifold with a finitely hyperbolic point. Intuition strongly suggests that a finitely hyperbolic point could only exist in a discrete space and not in a continuous space like an S-manifold. A peculiar property of lines in an S-manifold, however, is that a line that passes through a hyperbolic vertex is isolated from lines that are nearby (see Figure 2a). This ability to isolate lines suggests that it may be possible to construct an S-manifold with a finitely hyperbolic point.

**Acknowledgements**

Thanks to the Yahoo Smarandache Geometry Club for the interesting discussions and ideas that lead to this paper.

Thanks also to Ken Sullins for introducing me to ZAKS blocks and letting me play with his. Ohio Art's ZAKS blocks were very important in the development and presentation of these models.